\documentclass[12pt]{amsart}
\usepackage{amssymb}
\usepackage{amscd}
\usepackage{verbatim}
\usepackage{amssymb,amsfonts,pstricks,epsf}
\usepackage{graphicx}
\usepackage{epsfig}
\usepackage{color}
\usepackage{stmaryrd}

\newtheorem{theorem}{Theorem}[section]

\newtheorem{corollary}[theorem]{Corollary}
\newtheorem{innercustomthm}{Lemma}
\newenvironment{customthm}[1]
  {\renewcommand\theinnercustomthm{#1}\innercustomthm}
  {\endinnercustomthm}
  
\newtheorem{Innercustomthm}{Theorem}
\newenvironment{Customthm}[1]
  {\renewcommand\theInnercustomthm{#1}\Innercustomthm}
  {\endInnercustomthm}

\begin{document}

\title{Exit time moments and eigenvalue estimates}

\author{Emily B. Dryden, Jeffrey J. Langford \and Patrick McDonald}

\address{Department of Mathematics, Bucknell University, Lewisburg, Pennsylvania 17837}

\email{emily.dryden@Bucknell.edu}
\email{jeffrey.langford@Bucknell.edu}

\address{Division of Natural Science, New College of Florida, Sarasota, FL 34243}

\email{mcdonald@ncf.edu}

\date{\today}

\begin{abstract}

We give upper bounds on the principal Dirichlet eigenvalue associated to a smoothly bounded domain in a complete Riemannian manifold; the bounds involve $L^1$-norms of exit time moments of Brownian motion.  Our results generalize a classical inequality of P\'olya. We also prove lower bounds for Dirichlet eigenvalues using invariants that arise during the examination of the relationship between the heat content and exit time moments.\end{abstract}



\maketitle

\section{Introduction}

Literature devoted to the relationship between the geometry of a Riemannian manifold and the spectrum of the associated Laplace operator can be traced to Lord Rayleigh's conjecture for the principal eigenvalue associated to the Dirichlet Laplacian for planar domains.   Subsequent estimates for the principal Dirichlet eigenvalue in terms of geometric data associated to a given domain fill many journal pages.  The purpose of this paper is to develop a new family of such estimates that involve the relationship between the Dirichlet spectrum and the exit time moments of Brownian motion.  

We take as our point of departure the recent work of M.~van den Berg et al. \cite{BBV}, \cite{BFNT} which involves the relationship between the principal Dirichlet eigenvalue and {\it torsional rigidity}, an invariant that first arose in the nineteenth century study of elastica and that can be defined in terms of a solution of a Poisson problem.   More precisely, if $\Omega$ is a bounded open domain in a complete Riemannian manifold, then the torsional rigidity associated to $\Omega$ is 
\[
T_1(\Omega) = \int_\Omega u_1(x) dx,
\]
where $u_1$ solves the Poisson problem
\begin{eqnarray}
  \Delta u_1& = &-1 \hbox{ in } \Omega,  \label{T1.1}\\
 u_1 & = & 0 \hbox{ on } \partial \Omega. \label{T1.2}
\end{eqnarray}

For domains in Euclidean space, van den Berg, Buttazzo and Velichkov \cite{BBV} begin with P\'olya's inequality
\begin{equation}\label{inequality1}
\lambda_1(\Omega) T_1(\Omega) \leq |\Omega|,
\end{equation}
where $\lambda_1$ denotes the principal Dirichlet eigenvalue and $|\Omega |$ denotes volume.   They then study the domain functional 
\begin{equation}\label{domainfunctional}
F(\Omega) = \frac{\lambda_1(\Omega) T_1(\Omega)}{ |\Omega |}
\end{equation}
and note that the obvious upper bound is not sharp.  This work is continued in \cite{BFNT} where the authors improve the classical P\'olya bound with the estimate 
\begin{equation}\label{domainfunctional2}
F(\Omega) \leq 1- \nu_n T_1(\Omega) |\Omega |^{-1-\frac{2}{n}},
\end{equation}
where the constant $\nu_n$ depends only on dimension.  

Our work begins by noting that torsional rigidity can be expressed in terms of the exit time of Brownian motion.  More precisely, suppose $X_t$ is Brownian motion on a complete Riemannian manifold $M$ and $\Omega$ is an open domain in $M.$  Denote the first exit time of $X_t$ from $\Omega$ by $\tau:$
\[
\tau = \inf \{t\geq 0: X_t \notin \Omega\}.
\]
The torsional rigidity of $\Omega$ can then be expressed as 
\[
T_1(\Omega) = \int_\Omega {\mathbb E}^x [\tau] dx,
\]
where ${\mathbb E}^x $ denotes expectation with respect to the probability measure charging paths beginning at $x$ and $dx$ denotes the metric density.  Our main result involves extending the inequality (\ref{inequality1}) to higher moments and a more general context:  

\begin{theorem}\label{theorem1}  Let $M$ be a complete Riemannian manifold and let $\Omega \subset M$ be a bounded open domain.  Let $\tau$ be the first exit time of Brownian motion from $\Omega$ and let $T_k(\Omega)$ be the $L^1$-norm of the $k^{\mbox{th}}$ moment of $\tau:$   
\begin{equation}\label{moments1}
T_k(\Omega) = \int_\Omega {\mathbb E}^x [\tau^k] dx,
\end{equation}
where the integration is with respect to the measure induced by the Riemannian structure.  Then for every $k,$ 
\begin{equation}\label{inequalityk}
\lambda_1(\Omega)\leq \frac{(k!)^2}{(2k-1)!} \frac{T_{2k-1}(\Omega)}{T_{k}^2(\Omega)}|\Omega|.
\end{equation}
\end{theorem}

We call the sequence defined by (\ref{moments1}) the {\it $L^1$-moment spectrum associated to the domain $\Omega.$}  It is known that the $L^1$-moment spectrum determines the heat content asymptotics and, generically, the Dirichlet spectrum \cite{MM1}.  Recent work suggests the moment sequence determines a number of  geometric invariants of the associated domain and is of value in establishing comparison results (cf. \cite{HMP1}, \cite{HMP2}, \cite{HMP3}, \cite{Mc2}).  

The techniques required to obtain inequality (\ref{inequalityk}) are easily modified to obtain a variety of other inequalities.   For example, we prove:
 
 \begin{theorem}\label{theorem2}  Let $M$ and $\Omega$ be as in Theorem \ref{theorem1} and suppose $\{T_k(\Omega)\}$ is the $L^1$-moment spectrum associated to $\Omega.$  Then 
\begin{equation}\label{compk}
\lambda_1(\Omega) \leq  2k \frac{T_{2k-1}(\Omega)}{T_{2k}(\Omega)} .
\end{equation}
\end{theorem}

For the special case of geodesic balls in warped product spaces, Hurtado et al. \cite{HMP3} use (\ref{compk}) to prove a number of comparison results (for related results, see \cite{MM1}, \cite{Mc2}).  In addition, they establish that the asymptotics of the ratio given by the right hand side of inequality (\ref{compk}) determine the first Dirichlet eigenvalue.  There is an interpretation of this in terms of shape optimization and P\'olya's inequality in warped product spaces: the functional 
\[
F(\Omega) = \lim_{k\to \infty} \lambda_1(\Omega) \frac{T_{2k}(\Omega)}{2kT_{2k-1}(\Omega)}
\]
is optimized when $\Omega$ is a geodesic ball.

In addition to upper bounds for the principal Dirichlet eigenvalue, we study lower bounds for Dirichlet eigenvalues using invariants that arise during the examination of the relationship between the heat content and the moment spectrum.  More precisely, the $L^2$-norms of orthogonal projection of the constant function 1 on Dirichlet eigenspaces are closely related to the volume of the underlying domain, and there are natural lower bounds for eigenvalues in terms of these norms and the moment spectrum.  For example, we prove:

\begin{theorem}\label{theorem3} Let $M$ and $\Omega$ be as in Theorem \ref{theorem1}. Suppose that $\lambda$ is a Dirichlet eigenvalue associated to the domain $\Omega$ and let $E_\lambda$ be the eigenspace corresponding to $\lambda.$   Suppose that $a_\lambda^2$ is the square of the $L^2$-norm of the orthogonal projection of the constant function 1 on the eigenspace $E_\lambda.$   Then 
 \begin{equation}
 \left(\frac{k! a_{\lambda}^2  }{T_k(\Omega)}\right)^{\frac{1}{k}}  \leq  \lambda.
\end{equation}
\end{theorem}

Theorem \ref{theorem3} (see also Corollary \ref{lowerbound2}) strongly suggests that the invariants $a_\lambda^2$ play an important role in the geometric analysis of Riemannian spaces.  A review of the literature suggests that there is previous work, such as the Payne-Rayner inequality \cite{PR} and its subsequent refinements,  to support this view (see \cite{GN} for a survey).  

In addition to estimating eigenvalues using the $L^1$-norms of exit time moments, we investigate the relationship between eigenvalues and the $L^1$-norms of the variance of exit time (and the variance of higher powers of exit time).  Using the $L^1$-norm of variance we establish upper and lower bounds for Dirichlet eigenvalues (see Corollary \ref{varcor} and Corollary \ref{lowerbound2}).

The remainder of the paper is structured as follows.  In the next section we establish required notation and background information.  In Section \ref{Proofs} we complete the proofs of Theorems \ref{theorem1}-\ref{theorem3}, as well as related corollaries.  In the final section of the paper we discuss extensions of our results to operators which are of divergence form.     

\section{Background and notation}

As in the Introduction, let $M$ be a complete Riemannian manifold and suppose $X_t$ is Brownian motion on $M$ with infinitesimal generator $\Delta$.  Let ${\mathbb P}^x$ denote the probability measure charging Brownian paths beginning at $x$ and let ${\mathbb E}^x$ denote expectation with respect to ${\mathbb P}^x.$  For $\Omega\subset M$ a bounded open domain, let $\tau = \tau_{\Omega}$ be the first exit time of $X_t$ from $\Omega:$
\[
\tau = \inf\{t\geq 0: X_t \notin \Omega\}.
\]
There is a well-known connection between path properties of Brownian motion and solutions of the heat equation.  In particular, if 
\begin{equation}\label{heat1}
u(x,t)  = {\mathbb P}^x(\tau > t)
\end{equation}
then $u(x,t)$ satisfies
\begin{eqnarray}
\frac{\partial u}{\partial t} & = & \Delta u \hbox{ in } \Omega \times(0,\infty), \label{heat2}\\
u(x,0) & = & 1 \hbox{ in } \Omega, \label{heat3}\\
\lim_{x\to \sigma} u(x,t) &= & 0 \hbox{ for all }(\sigma,t)  \in \partial \Omega \times (0,\infty) \label{heat4}.
\end{eqnarray}
The moments of the exit time are given by integration against the appropriate probability distribution, which is given by the solution of (\ref{heat2})-(\ref{heat4}).  More precisely, since $u(x,t) = {\mathbb E}^x[1_{\tau^k> t^k}],$ we have 
\begin{equation}
{\mathbb E}^x [\tau^k] = k\int_0^\infty t^{k-1} u(x,t) dt.
\end{equation}
If we integrate over the domain and use Fubini, we obtain 
\begin{equation}\label{momentsheat}
T_k(\Omega) = k\int_0^\infty t^{k-1} H(t) dt,
\end{equation}
where $T_k(\Omega)$ is the $k^{\textup{th}}$ element of the moment spectrum and $H(t)$ is the {\it heat content} of the domain:
\begin{equation}
H(t) = \int_\Omega u(x,t) dx.
\end{equation}

Expressing the heat content in terms of spectral data provides a relationship between Dirichlet spectrum and the $L^1$-moments of the exit time of Brownian motion.  More precisely, given the Dirichlet spectrum associated to $\Omega,$ fix an orthonormal basis of corresponding eigenfunctions.  Then the heat content associated to $\Omega$ can be written as 
\begin{equation}
H(t) = \sum_{\lambda\in \hbox{spec}^*(\Omega)} a^2_\lambda e^{-\lambda t},
\end{equation}
where $\hbox{spec}^*(\Omega)$ is the set of Dirichlet eigenvalues (omitting multiplicity) and $a_\lambda^2$ is the square of  the $L^2$-norm of the orthogonal projection of the constant function $1$ on the eigenspace corresponding to eigenvalue $\lambda.$   Summing over those $\lambda$ for which the corresponding eigenspace is not orthogonal to constant functions, we see that the  constants $a_\lambda^2$ are invariants of the isometry group which are closely related to the volume.  In fact, using the asymptotic behavior of the heat content for small time, we have 
\begin{equation}\label{volumepartition}
|\Omega | \sim_{t \to 0^+}  \sum_{\lambda\in \hbox{spec}^*(\Omega)} a^2_\lambda
\end{equation}
and the $a_\lambda^2$ partition the volume of the domain $\Omega$ among the eigenspaces.  

Equation (\ref{momentsheat}) should be interpreted as connecting the $L^1$-moment spectrum to the Mellin transform of the heat content.   As noted in \cite{MM1}, the Mellin transform of the heat content takes the form of a Dirichlet series
\begin{equation}
\zeta_{\Omega} (s) =  \sum_{\lambda\in \hbox{spec}^*(\Omega)} a^2_\lambda \left( \frac{1}{\lambda }\right)^s
\end{equation}
and the $L^1$-moment spectrum is related to the Dirichlet spectrum via the equality
\begin{equation}\label{zeta}
\Gamma(k+1) \zeta_\Omega(k) = T_k(\Omega),
\end{equation}
 where $\Gamma(k)$ is the gamma function.  There is a great deal of information contained in the relationship (\ref{zeta}).  For example, (\ref{zeta}) implies that the asymptotics of the moment spectrum determine the principal Dirichlet eigenvalue (cf. \cite{MM1}):
 \begin{equation}\label{MM11}
 \lambda_1 = \sup\left\{ \eta: \limsup_{n\to \infty} \left(\eta^n\frac{ T_n(\Omega)}{\Gamma(n+1)}\right) < \infty \right\}.
 \end{equation}
 Because (\ref{MM11}) indicates that the tail of the $L^1$-moment spectrum determines $\lambda_1,$ it also determines $a_{\lambda_1}^2:$
 \begin{equation}\label{a1L}
 a_{\lambda_1}^2 = \limsup_{n\to \infty} \left\{\lambda_1^n \frac{ T_n(\Omega)}{\Gamma(n+1)} \right\}.
 \end{equation}
 
As a final bit of background information, we note that the moment spectrum can be written in terms of a hierarchy of Poisson problems, a fact that will be very useful in obtaining relationships among the moment spectrum, the Dirichlet spectrum and the geometry of the domain.  In more detail, suppose $u_1$ is defined as the solution of the Poisson problem given in (\ref{T1.1})-(\ref{T1.2}).  Define $u_k$ inductively as the solution of the Poisson problem
\begin{eqnarray}
\Delta u_k& = & -ku_{k-1} \hbox{ in } \Omega \label{recursion11}, \\
u_k & = & 0 \hbox{ on } \partial \Omega.\label{recursion2}
\end{eqnarray}

Then $u_k(x) = {\mathbb E}^x[\tau^k]$ (cf. \cite{Mc1}, \cite{KMM} for details).    

\section{Proofs}\label{Proofs}
We begin this section with a combined proof of Theorems \ref{theorem1} and \ref{theorem2}:

\begin{proof}
Recall that the first Dirichlet eigenvalue minimizes the Rayleigh quotient:
\[
\lambda_1(\Omega) = \inf_{\phi \in H^1_0(\Omega)} \frac{ \int_\Omega |\nabla \phi|^2 dx }{ \int_\Omega \phi^2 dx }.
\]
We take $\phi = u_k$ as a test function and note
\begin{eqnarray*}
\int_\Omega |\nabla u_k|^2 dx & = & -\int_\Omega u_k \Delta u_k dx \\
 & = & k \int_\Omega u_k  u_{k-1} dx \\
 & = & -\frac{k}{k+1} \int_\Omega \Delta u_{k+1}  u_{k-1} dx.
 \end{eqnarray*}
Iterating this procedure and using the self-adjointness of $\Delta$, we obtain

 \begin{equation}\label{th1e1}
 \lambda_1(\Omega)  \leq \frac{(k!)^2}{(2k-1)!}\frac{ \int_\Omega  u_{2k-1}  dx}{\int_\Omega u_k^2 dx}.
\end{equation}
By H\"older's inequality, 
\begin{equation} \label{th1e3}
\left(\int_\Omega u_k dx\right)^2   \leq |\Omega| \int_\Omega u_k^2 dx .
\end{equation}
Using (\ref{th1e3}) in (\ref{th1e1}) completes the proof of Theorem \ref{theorem1}.

To prove Theorem \ref{theorem2}, we return to the denominator in (\ref{th1e1}), but instead of using H\"older's inequality we integrate by parts.  We begin by noting 
\begin{eqnarray*}
\int_\Omega u_k u_k dx & = & -\frac{1}{k+1} \int_\Omega \Delta u_{k+1} u_k dx\\
 & = & \frac{k}{k+1} \int_\Omega u_{k+1} u_{k-1} dx .
 \end{eqnarray*}
Iterating, we obtain
\begin{equation}\label{variance2}
\int_\Omega u_k ^2 dx  = \frac{(k!)^2}{(2k)!} \int_\Omega u_{2k} dx.
\end{equation}
Substituting (\ref{variance2}) in (\ref{th1e1}) completes the proof of Theorem \ref{theorem2}.
\end{proof}

The ratio of integrals appearing on the right hand side of inequality (\ref{th1e1}) has an interesting interpretation in terms of Brownian motion.  To make this clear, we compute the variance of the random variable $\tau^k:$
\begin{eqnarray}
\hbox{Var}[\tau^k] &=& {\mathbb E}^x[(\tau^k -{\mathbb E}^x[\tau^k])^2] \nonumber \\
 & = &  {\mathbb E}^x[\tau^{2k} -u_k^2] \nonumber \\
  & = & u_{2k}(x) - u_k^2(x).\label{variance1}
  \end{eqnarray}
Integrating the variance over the domain produces an $L^1$-norm of the variance of $\tau^k,$ which we denote by $\hbox{Var}_k(\Omega):$
\begin{equation}\label{vardef}
\hbox{Var}_k(\Omega)  = \int_\Omega (u_{2k}- u_k^2)  dx.
\end{equation}
Using (\ref{variance2}) we have 
\[
\int_\Omega u_{2k} dx = \frac{(2k)!}{(k!)^2} \int_\Omega u_k^2 dx.
\]
We use this to rewrite the $L^1$-norm of the variance of $\tau^k$ in terms of $\int_\Omega u_k^2 dx:$ 
\begin{equation}\label{variance4}
\int_\Omega u_k^2 dx = \frac{(k!)^2}{(2k)! - (k!)^2} \hbox{Var}_k(\Omega).
\end{equation}
Substituting (\ref{variance4}) in (\ref{th1e1}) we conclude:

\begin{corollary}\label{varcor} Let $\Omega$ be as in Theorem \ref{theorem1}. For $k$ a positive integer, let $\textup{Var}_k(\Omega)$ be the $L^1$-norm of the variance of $\tau^k:$ 
\[
\textup{Var}_k(\Omega) = \int_\Omega (u_{2k} - u_k^2) dx.
\] 
Then
\begin{equation}\label{varestimate}
\lambda_1(\Omega) \leq  \frac{(2k)! - ( k!)^2 }{(2k-1)!} \frac{T_{2k-1}(\Omega)}{\textup{Var}_k(\Omega)}.
\end{equation}
\end{corollary}

Next, we can establish the lower bounds of Theorem \ref{theorem3}:

\begin{proof}
To obtain lower bounds we use a variational formula for $T_k(\Omega)$ (see the appendix for details):
\[
T_k(\Omega)= 
\begin{cases}
k!\underset{\phi \in \mathcal F_k}{\sup} \frac{\left(\int_{\Omega}\phi dx\right)^2}{\int_{\Omega} (\Delta^{\frac{k}{2}}\phi)^2dx} & \textup{ if $k$ is even},\\
k!\underset{\phi \in \mathcal F_k}{\sup} \frac{\left(\int_{\Omega}\phi dx\right)^2}{\int_{\Omega} |\nabla \Delta^{\llbracket \frac{k}{2}\rrbracket}\phi|^2dx} & \textup{ if $k$ is odd}.
\end{cases}
\]
Here $\llbracket \cdot \rrbracket$ denotes the greatest integer function and $\mathcal F_k$ denotes the class of functions
\[
\mathcal F_k=\{\phi \in C^{\infty}(\overline \Omega):\int_{\Omega}\phi dx\neq 0,\  \Delta^j \phi=0 \textup{ on }\partial \Omega \textup{ for }0\leq j \leq k-1\}.
\]
Let $\phi$ be the $L^2$-normalized first Dirichlet eigenfunction which is positive.  Taking $\phi$  as a test function in the variational quotient for $T_k(\Omega)$, we obtain in both cases:
\begin{equation}\label{lower1}
\lambda_1^k T_k(\Omega) \geq k! a_{\lambda_1}^2  ,
\end{equation}
where
\[
a_{\lambda_1}^2 = \left(\int_\Omega \phi dx\right)^2 
\]
is the square of the $L^2$-norm of the projection of the constant function 1 on the eigenspace corresponding to the principal Dirichlet eigenvalue.  This proves Theorem \ref{theorem3} for the principal Dirichlet eigenvalue.     

Because the principal Dirichlet eigenvalue is simple and the corresponding eigenfunction has a sign, the square of the $L^2$-norm of orthogonal projection of the function 1 coincides with the square of the $L^1$-norm of the $L^2$-normalized eigenfunction.  For higher eigenvalues this will not be the case, but there is a closely related estimate.  More precisely, suppose $\lambda $ is any Dirichlet eigenvalue with corresponding eigenspace $E_\lambda$ of dimension $n.$  Fix an orthonormal basis $\{\phi_i\}_{i=1}^n $ of $E_\lambda.$   Let $c_i= \int_\Omega \phi_i dx$ and let $\phi = \sum_{i=1}^n c_i \phi_i$ be orthogonal projection of the constant function 1 on the eigenspace $E_\lambda.$  Then using orthogonality,
\begin{equation}\label{L2prop}
\int_\Omega \phi dx = \int_\Omega \phi^2 dx.
\end{equation}
As above, we write the square of the $L^2$-norm of $\phi$ as $a_\lambda^2.$  Then, assuming $E_\lambda$ is not orthogonal to constants, the variational quotient gives 
\begin{equation}\label{lower2}
\lambda^k T_k(\Omega) \geq k! a_{\lambda}^2.
\end{equation}
Since inequality (\ref{lower2}) clearly holds if $E_\lambda$ is orthogonal to constants, Theorem \ref{theorem3} follows.
\end{proof}

As was the case with upper bounds, there are lower bounds of Dirichlet eigenvalues involving the $L^1$-norm of the variance.  In fact, suppose $\lambda$ is as above and $\phi$ is the orthogonal projection of the constant function 1 on the eigenspace $E_\lambda$ corresponding to $\lambda.$  Arguing inductively,
\begin{eqnarray}
a_\lambda^2 & = & \int_\Omega \phi dx \nonumber \\
  & = & - \int_\Omega \phi \Delta u_1 dx \nonumber \\
  & = & \lambda \int_\Omega \phi u_1 dx \nonumber\\
  & = & \frac{\lambda^k}{k!} \int_\Omega \phi u_k dx. \label{lower3}
\end{eqnarray}
Applying H\"older, we obtain 
\begin{equation}\label{lower4}
\int_\Omega \phi u_k dx \leq  \left( \int_\Omega \phi^2 dx\right)^{\frac12}  \left( \int_\Omega u_k^2 dx \right)^{\frac12} .
\end{equation}
Combining (\ref{lower4}), (\ref{lower3}) and (\ref{variance4}) we have 
\[
a_\lambda^2  \leq \frac{\lambda^k}{k!}a_\lambda \left (\frac{(k!)^2}{(2k)!-(k!)^2}\textup{Var}_k(\Omega)\right)^{\frac12}.
\]
This proves:
\begin{corollary}\label{lowerbound2} Let $\Omega$ be as in Theorem \ref{theorem1} and suppose $\lambda $ is a Dirichlet eigenvalue with associated eigenspace not orthogonal to constants.  Then
\begin{equation}
\left (a_\lambda^2\frac{(2k)!-(k!)^2}{\textup{Var}_k(\Omega)} \right)^{\frac{1}{2k}} \leq \lambda(\Omega),
\end{equation}
where $a_\lambda^2$ is the square of the $L^2$-norm of the orthogonal projection of the constant function 1 on the eigenspace corresponding to $\lambda$ and $\textup{Var}_k(\Omega)$ is as in Corollary \ref{varcor}.
\end{corollary}

\section{Divergence form operators}

Let $M$ be a complete Riemannian manifold and denote by $\nabla$ the gradient operator.  A second-order partial differential operator $L$ acting on smooth functions on $M$ is said to be of divergence form if there is a symmetric endomorphism of the tangent bundle, $A,$ such that the action of $L$ is given by 
\[
Lf = \hbox{div}(A\nabla f),
\]
where $\hbox{div}$ is the divergence operator acting on vector fields.   We will assume throughout this section that $L$ is of divergence form and uniformly elliptic (in particular, the spectrum of $A$ is bounded below by a positive constant).  

Let $\Omega\subset M$ be a smoothly bounded domain with compact closure, and $X$ be a smooth vector field on $\Omega$.  Then the Divergence Theorem holds:
\[
\int_\Omega \hbox{div}(X) dx = \int_{\partial \Omega} \langle X, n \rangle dS,
\]
where $dx$ is the metric density, the pairing is given by the metric and $dS$ is the induced measure on the boundary.  In particular, if we choose $f$ to vanish on the boundary of $\Omega,$ we have 
\[
\int_\Omega L(f^2) dx = 0
\]
from which we obtain 
\begin{equation}
\int_\Omega |\nabla f|_A^2 dx = \int_\Omega fLf dx,
\end{equation}
where the norm $|\nabla f|_A^2$ is defined via pairing with the metric:
\[
|\nabla f|_A^2 = -\langle A\nabla f,\nabla f \rangle.
\]
It follows that the operator $L$ acting on smooth functions that vanish at the boundary of a smoothly bounded domain with compact closure is formally self-adjoint (that is, we can integrate by parts).  In particular, the standard functional analytic tools are available, and we conclude that, as in the case of the Dirichlet Laplacian, the spectrum of $L$ will be positive, discrete and accessible via Rayleigh quotients.   

If $L$ is as above, we can associate to $L$ a hierarchy of Poisson problems in analogy to (\ref{T1.1})-(\ref{T1.2}) and (\ref{recursion11})-(\ref{recursion2}):  Suppose $u_1$ solves 
\begin{eqnarray}
Lu_1 & = & -1 \hbox{ in } \Omega, \label{L11}\\
u_1 & = & 0 \hbox{ on } \partial \Omega, \label{L12}
\end{eqnarray}
and having defined $u_{k-1},$ define $u_k$ by  
\begin{eqnarray}
Lu_k & = & -ku_{k-1} \hbox{ in } \Omega, \label{L13}\\
u_k & = & 0 \hbox{ on } \partial \Omega.\label{L14}
\end{eqnarray}

The relationships established in Theorem \ref{theorem1} and Theorem \ref{theorem2} depend only upon integration by parts, H\"older's inequality, the functional analytic tools sketched above, and the definition of the Poisson hierarchy.  Thus, the proofs of these theorems carry over, mutatis mutandi, to establish the analogous theorems for divergence form operators.  More precisely, we have:

\begin{theorem}Suppose $M$ is a complete Riemannian manifold and $\Omega \subset M$ is a smoothly bounded domain with compact closure.  Suppose $L$ is a uniformly elliptic divergence form operator and let 
\[
T_{k,L}(\Omega) = \int_\Omega u_k dx,
\]
where $dx$ is the metric density and $u_k$ is defined by (\ref{L11})-(\ref{L14}).  Let $\lambda_{1,L}$ be the principal Dirichlet eigenvalue for $L.$  Then
\begin{equation}
\lambda_{1,L}(\Omega) \leq \frac{(k!)^2}{(2k-1)!}\frac{T_{2k-1,L}(\Omega)}{T_{k,L}^2(\Omega)}|\Omega |.
\end{equation}
In addition,
\begin{equation}
\lambda_{1,L}(\Omega) \leq 2k \frac{T_{2k-1,L}(\Omega)}{T_{2k,L}(\Omega)}.
\end{equation}
\end{theorem}

A careful review of the proofs of Theorem \ref{theorem3} and the corollaries of section 3 indicates that what is formally required is no more than machinery discussed above.  On the other hand, in order to give the references to ``exit time'' meaning, some work must be done.  In particular, we must establish a relationship between the Poisson hierarchy, the heat operator associated to $L$ and a description of the solution of the heat operator using path properties of some process.    When $L$ is the Laplace operator, the connection between the heat operator and the Poisson hierarchy hinges on the fact that the Laplace operator is the infinitesimal generator of the Brownian process.  In the present circumstances, there is analogous structure: there is a diffusion process whose infinitesimal generator is $L$ and for which the exit time moments from the domain $\Omega$ are given by the solution of the appropriate member of the Poisson hierarchy.  The details are discussed in \cite{Mc1} and references therein.   This being the case, the proofs of the results presented above carry over, mutatis mutandi, to the case of divergence form operators and the objects described can be discussed probabilistically.  While the squared norms of projections on eigenspaces will play the same algebraic and analytic role, there will, in general, be no relation to the underlying geometry of the domain $\Omega.$

\section*{Appendix}
Here we establish the variational characterization for $T_k(\Omega)$ used in the proof of Theorem \ref{theorem3}.  A similar characterization in the more general case of an elliptic operator of divergence form may be found in \cite{KM}. In what follows, we denote
\[
\mathcal F_k=\{\phi \in C^{\infty}(\overline \Omega):\int_{\Omega}\phi dx\neq 0,\  \Delta^j \phi=0 \textup{ on }\partial \Omega \textup{ for }0\leq j \leq k-1\}.
\]
We begin with a preliminary lemma:

\begin{customthm}{A1}\label{Var Lemma} Let $\Omega$ be as in Theorem \ref{theorem1} and $u_k$ be as in \eqref{recursion11} and \eqref{recursion2}. If $f\in \mathcal F_k$, then
\begin{equation}\label{Eqn:LemmaEven}
\int_{\Omega}fdx=\frac{(-1)^k}{k!}\int_{\Omega}(\Delta^kf)u_kdx.
\end{equation}
If $f\in \mathcal F_{k+1}$, then
\begin{equation}\label{Eqn:LemmaOdd}
\int_{\Omega}fdx=\frac{(-1)^k}{(k+1)!}\int_{\Omega}\nabla (\Delta^kf)\cdot \nabla u_{k+1}dx.
\end{equation}
\end{customthm}

\begin{proof}
To establish \eqref{Eqn:LemmaEven}, we directly compute
\begin{align*}
\int_{\Omega}fdx&=-\int_{\Omega}f\Delta u_1dx\\
&=-\int_{\Omega}(\Delta f)u_1dx\\
&=\frac{1}{2}\int_{\Omega}(\Delta f)(\Delta u_2)dx\\
&=\frac{1}{2}\int_{\Omega}(\Delta^2f)u_2dx \qquad \textup{(if $\Delta f=0$ on $\partial \Omega$)}.
\end{align*}
Iterating this argument establishes \eqref{Eqn:LemmaEven}. To prove \eqref{Eqn:LemmaOdd}, we use \eqref{Eqn:LemmaEven}:
\begin{align*}
\int_{\Omega}fdx&=\frac{(-1)^{k+1}}{(k+1)!}\int_{\Omega}u_{k+1} \Delta(\Delta^kf)dx\\
&=\frac{(-1)^k}{(k+1)!}\int_{\Omega}\nabla(\Delta^kf)\cdot \nabla u_{k+1}dx,
\end{align*}
where the last equality follows from a Green's identity and the fact that $\Delta^kf=0$ on $\partial\Omega$.
\end{proof}

We may now establish the desired variational characterization.
\begin{Customthm}{A2}
Let $\Omega$ be as in Theorem \ref{theorem1}. Then 
\[
T_k(\Omega)= 
\begin{cases}
k!\underset{\phi \in \mathcal F_k}{\sup} \frac{\left(\int_{\Omega}\phi dx\right)^2}{\int_{\Omega} (\Delta^{\frac{k}{2}}\phi)^2dx} & \textup{ if $k$ is even},\\
k!\underset{\phi \in \mathcal F_k}{\sup} \frac{\left(\int_{\Omega}\phi dx\right)^2}{\int_{\Omega} |\nabla \Delta^{\llbracket \frac{k}{2}\rrbracket}\phi|^2dx} & \textup{ if $k$ is odd},
\end{cases}
\]
where $\llbracket \cdot \rrbracket$ denotes the greatest integer function.
\end{Customthm}

\begin{proof}
For the even case, let $\mathcal G=\{\Delta^nf:f\in \mathcal F_{2n}\} \subseteq \mathcal F_n$ and define $Q:\mathcal G \to \mathbb{R}$ via
\begin{align*}
Q(\Delta^nf)&=\left(\frac{1}{n!}\right)^2\frac{\left(\int_{\Omega}u_n(\Delta^nf)dx\right)^2}{\int_{\Omega}(\Delta^nf)^2dx}\\
&=\left(\frac{1}{n!}\right)^2\frac{\langle u_n, \Delta^nf\rangle^2}{\|\Delta^nf\|_2^2},
\end{align*}
where $\langle\cdot,\cdot \rangle$ denotes the standard inner product on $L^2(\Omega)$ and $\|\cdot\|_2$ denotes $L^2$-norm. We observe that $Q$ is largest when $\Delta^nf$ is a scalar multiple of $u_n$. It follows that $Q$ is maximized when $f=u_{2n}$. Thus,
\begin{align*}
\sup_{f\in \mathcal F_{2n}}Q(\Delta^nf)&=\left(\frac{1}{n!}\right)^2\frac{\left(\int_{\Omega}u_n\Delta^nu_{2n}dx\right)^2}{\int_{\Omega}\left(\Delta^nu_{2n}\right)^2dx}\\
&=\frac{n!}{(-1)^n(2n)!}\frac{\left(\int_{\Omega}u_{2n}dx\right)^2}{\int_{\Omega}u_n\Delta^nu_{2n}dx}\\
&=\frac{1}{(2n)!}T_{2n}(\Omega),
\end{align*}
where the second and third equalities invoke \eqref{Eqn:LemmaEven}. Using \eqref{Eqn:LemmaEven} again, we see that 
\[
Q(\Delta^nf) = \frac{\left( \int_{\Omega} f dx \right)^2}{\int_{\Omega}(\Delta^n f)^2 dx},
\]
establishing the even case.

For the odd case, set $\mathcal G=\{\nabla (\Delta^nf):f\in \mathcal F_{2n+1}\}$ and define $Q:\mathcal G \to \mathbb{R}$ via
\[
Q(\nabla(\Delta^nf))=\left(\frac{1}{(n+1)!}\right)^2 \frac{\left(\int_{\Omega}\nabla u_{n+1}\cdot \nabla(\Delta^n f)dx\right)^2}{\int_{\Omega}\left(\nabla(\Delta^nf)\right)^2dx}.
\]
As before, $Q$ is maximized when $\nabla(\Delta^nf)$ is a scalar multiple of $\nabla u_{n+1}$. Taking $f=u_{2n+1}$ and invoking \eqref{Eqn:LemmaOdd}, the remainder of the proof proceeds in a manner parallel to the even case.
\end{proof}

\bibliographystyle{plain}
\bibliography{highermoments}

\end{document}